\documentclass[12pt]{article}
\usepackage{amssymb}
\usepackage{graphicx}
\usepackage{amsfonts}
\usepackage{latexsym}
\usepackage{amsthm}
\usepackage{amsmath}

\usepackage{amssymb, amscd}

\usepackage{url}
\usepackage[T1]{fontenc}

\newtheorem{problem}{Problem}
\newtheorem{theo}[problem]{Theorem}

\newtheorem{prob}[problem]{Problem}

\newtheorem{prop}[problem]{Proposition}
\newtheorem{cor}[problem]{Corollary}
\newtheorem{lema}[problem]{Lemma}

\newtheorem{conj}[problem]{Conjecture}

\begin{document}
\date{}
 \title{{Topological obstructions to totally skew embeddings}}

\author{{\DJ{}or\dj{}e Barali\' c}\\ {\small Mathematical Institute}\\[-2mm] {\small SASA,
Belgrade} \and Branislav Prvulovi\' c\\ {\small Faculty of Mathematics}\\[-2mm] {\small University of Belgrade}
\and Gordana Stojanovi\' c\\ {\small Faculty of Mathematics}\\[-2mm] {\small University of Belgrade}
\and Sini\v sa Vre\' cica\\ {\small Faculty of Mathematics}\\[-2mm] {\small University of Belgrade}
  \and Rade  \v Zivaljevi\' c\\ {\small Mathematical Institute}\\[-2mm] {\small SASA,
  Belgrade}\\[-2mm]}

\maketitle
\begin{abstract}
Following Ghomi and Tabachnikov \cite{Gho-Tab} we study the
invariant $N(M^n)$ defined as the smallest dimension $N$ such that
there exists a {\em totally skew} embedding of a smooth manifold
$M^n$ in $\mathbb{R}^N$. This problem is naturally  related to the
question of estimating the geometric dimension of the stable
normal bundle of the configuration space $F_2(M^n)$ of ordered
pairs of distinct points in $M^n$. We demonstrate that in a number
of interesting cases the lower bounds on $N(M^n)$ obtained by this
method are quite accurate and very close to the best known general
upper bound $N(M^n)\leq 4n+1$ established in \cite{Gho-Tab}. We
also provide some evidence for the conjecture that for every
$n$-dimensional, compact smooth manifold $M^n$ $(n>1)$,
$$N(M^n)\leq 4n-2\alpha (n)+1.$$
\end{abstract}

\renewcommand{\thefootnote}{}
\footnotetext{This research was supported by the Grants 144014 and
144026 of the Ministry for Science and Technological Development
of Serbia.}

\section{Introduction}
\label{sec:intro}

Two lines in an affine space $\mathbb{R}^N$ are called {\em skew}
if they are neither parallel nor have a point in common or
equivalently if their affine span has dimension $3$.  More
generally, affine subspaces $U_1,\ldots, U_l$ of $\mathbb{R}^N$
are called {\em skew} if their affine span has dimension ${\rm
dim}(U_1)+\cdots +{\rm dim}(U_l)+l-1$, in particular a pair $U,V$
of affine subspaces of $\mathbb{R}^N$ is skew if and only if each
two lines $p\subset U$ and $q\subset V$ are skew.

\medskip
An embedding $f : M^n\rightarrow \mathbb{R}^{N}$ of a smooth
manifold is called {\em totally skew} if for each two distinct
points $x,y\in M^n$ the affine subspaces $df(T_xM)$ and $df(T_yM)$
of\, $\mathbb{R}^N$ are skew. Define $N(M^n)$ as the minimum $N$
such that there exists a totally skew embedding of $M^n$ into
$\mathbb{R}^N$.

\medskip
Ghomi and Tabachnikov began in \cite{Gho-Tab} the study of totally
skew embeddings of mani\-folds and established a surprising
connection of $N(M^n)$ with some classical invariants of geometry
and topology. For example they showed \cite[Theorem~1.4]{Gho-Tab}
that the problem of estimating $N(\mathbb{R}^n)$ is intimately
related to the generalized vector field problem and the immersion
problem for real projective spaces,  as exemplified by the
inequality
$$
N(\mathbb{R}^n)\geq r(n) +n
$$
where $r(n)$ is the minimum $r$ such that the Whitney sum
$r\xi_{n-1}$ of $r$ copies of the canonical line bundle over
$\mathbb{R}P^{n-1}$ admits $n+1$ linearly independent continuous
cross-sections.

\medskip
Another example (\cite[Theorem~1.2]{Gho-Tab}) is the inequality
$$
N(S^n)\leq n + m(n) +1
$$
where $m(n)$ is an equally well-known function defined as the
minimum $m$ such that there exists a non-singular, symmetric
bilinear form $B : \mathbb{R}^{n+1}\times
\mathbb{R}^{n+1}\rightarrow \mathbb{R}^m$. As a consequence they
deduced the inequalities $N(S^n)\leq 3n+2$ and $N(S^{2k+1})\leq
3(2k+1)+1$.

\medskip
It appears that very little is known about the exact values of
$N(M)$. Indeed, according to \cite{Gho-Tab}, the only currently
known exact values of this invariant are,
$$N(\mathbb{R}^1) = 3,\quad N(S^1)=4,\quad N(\mathbb{R}^2)=6.$$
Finally for a general $n$-manifold $M^n$ Ghomi and Tabachnikov
established upper and lower bounds
\begin{equation}\label{eqn:lower-upper}
2n+1\leq N(M^n)\leq 4n+1
\end{equation}
and showed that the lower bound can be improved to $2n+2$ if $M^n$
is a closed manifold.

\bigskip
In this paper we are interested in topological obstructions to
totally skew embeddings of manifolds, in particular we address the
problem of finding good lower bounds for $N(M^n)$. We demonstrate
that in many classes of manifolds there are examples where the
upper bound $4n+1$ from (\ref{eqn:lower-upper}) is very close to
the actual value of $N(M^n)$. For example $N(\mathbb{R}{P}^n)$ is
by Proposition~\ref{prop:proj-donja-ocena} one of the numbers
$4n-1, 4n, 4n+1$ if $n=2^k$ is a power of $2$, in particular
$N(\mathbb{R}P^2)$ is 7, 8, or 9. More generally, if $M^n =
\mathbb{R}P^{n_1}\times\cdots\times \mathbb{R}P^{n_k}$ is a
product of real projective spaces and $n_i=2^{r_i}$ are different
powers of $2$, then (Theorem~\ref{thm:proj-product})
$$
N(M^n) = N(\mathbb{R}{P}^{n_1}\times\cdots\times
\mathbb{R}{P}^{n_k}) \geq 4n - 2\alpha(n) + 1
$$
where $\alpha(n)$ is number of non-zero digits in the binary
representation of $n$. A similar bound
(Theorem~\ref{thm:compl-proj-product})
$$
N(X)\geq 8n-4\alpha (n)+1 = 4\cdot{\rm dim}_{\mathbb{R}}(X) -
4\alpha(n) +1
$$
is obtained if $X=\mathbb{C}P^{n_1}\times \cdots \times
\mathbb{C}P^{n_k}$ where $n_i=2^{r_i}$ are different powers of $2$
and $n=n_1+\cdots +n_k={\rm dim}_{\mathbb{C}}(X)$.

\medskip In pursuit of other examples of manifolds where $N(M^n)$
gets very close to the upper bound $4n+1$ we continue with the
analysis of Grassmann manifolds $G_k(\mathbb{R}^{n+k})$ and their
oriented counterparts $\tilde{G}_k(\mathbb{R}^{n+k})$. For example
(Theorems~\ref{t1} and \ref{t3}) we prove that
$N\left(G_2\left(\mathbb{R}^{2^r+2}\right)\right)\geq 4\cdot
2^{r+1}-3$ and $N(\tilde{G}_2(\mathbb{R}^{2^r+2}))\geq 3\cdot
2^{r+1}+1$. Similar inequalities can be expected for many other
Grassmannians as illustrated by the inequalities
$$
N(G_3(\mathbb{R}^6))\geq 31,\quad N(G_3(\mathbb{R}^7))\geq
43,\quad N(\tilde{G}_3(\mathbb{R}^7))\geq 41, \mbox{ {\rm etc.} }
$$
These results are in sharp contrast with the fact that very little
is known about the exact values of the invariant $N(M^n)$, for
example the exact value of $N(M^2)$ is not known for any closed
surface $M^2$. This and a sample of other open problems and
conjectures can be found in the final section of the paper where
we also offer a brief outlook to future research.

Possibly the most intriguing and attractive is the
Conjecture~\ref{conj:najinteresatnije} which, in analogy with the
classical {\em Immersion Conjecture} \cite{Cohen}, predicts that
for $n>1$
$$
N(M^n)\leq 4n - 2\alpha(n)+1.
$$

\section{Vector bundle decomposition}
\label{sec:vector-bundle}

Let $F_2(M):= M^2\setminus\Delta_M$ be the configuration space
(manifold) of all distinct ordered pairs of points in $M$. The
tangent bundle $T(F_2(M))$ admits a splitting
\begin{equation}\label{eqn:splitting}
T(F_2(M))\cong \pi_1^\ast TM\oplus \pi_2^\ast TM
\end{equation}
where $\pi_1, \pi_2 : F_2(M)\rightarrow M$ are the natural
projections. Simplifying the notation let $T_{(x,y)}(F_2(M)) \cong
T_x(M)\oplus T_y(M)$ be the fibre of this bundle at $(x,y)\in
F_2(M)$.

\begin{figure}[hbt]
\centering
\includegraphics[scale=0.45]{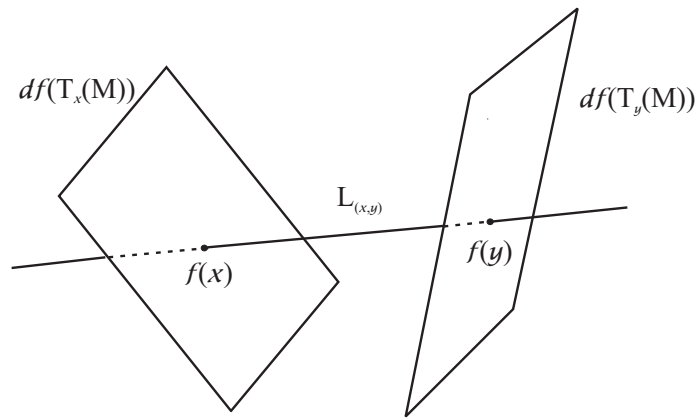}
\caption{Fibre of the bundle
$T(F_2(M))\oplus\varepsilon^1$.}\label{fig:prvi-kvadrant}
\end{figure}

\medskip
If $f : M^n\rightarrow \mathbb{R}^{N}$ is an embedding, then there
is a trivial line bundle $L$ over $F_2(M)$ such that for $(x,y)\in
F_2(M)$ the fibre $L_{(x,y)}$ is the line $\mathbb{R}\cdot
(f(y)-f(x))$. If $f : M^n\rightarrow \mathbb{R}^{N}$ is a totally
skew embedding, then there arises a monomorphism of vector bundles
$$
\Phi = \Phi^{(1)}\oplus \Phi^{(2)}: T(F_2(M))\oplus \varepsilon^1
\longrightarrow F_2(M)\times \mathbb{R}^N
$$
where $\Phi^{(1)}_{(x,y)} : T_x(M)\oplus T_y(M) \rightarrow
\mathbb{R}^N$ is the map defined by $\Phi^{(1)}_{(x,y)}(u,v) =
df_x(u) + df_y(v)$ and $\Phi^{(2)}$, defined by
$\Phi^{(2)}(\lambda)=\lambda(f(y)-f(x))$, maps the trivial line
bundle $\varepsilon^1$ to $L$. In this case the trivial
$N$-dimensional bundle $\varepsilon^N$ over $F_2(M)$ splits as
 \begin{equation}\label{eqn:splits}
\varepsilon^N \cong T(F_2(M))\oplus \varepsilon^1\oplus \nu
 \end{equation}
where $\nu$ is a $(N-2n-1)$-dimensional ``normal'' bundle. As a
consequence (\cite[Section~4]{Mi-Sta}) we obtain the following
proposition.

\begin{prop}\label{prop:dual}
If the dual Stiefel-Whitney class
$$\overline{w}_k(T(F_2(M))) := w_k(\nu)\in
H^k(F_2(M),\mathbb{F}_2)$$ is non-zero, then $2n+k+1\leq N$. In
particular, $N(M)\geq 2n+k+1$.
\end{prop}

\section{Characteristic classes of $T(F_2(M))$}
\label{sec:characteristic-classes}

The cohomology of $F_2(M)= M^2\setminus \Delta_M$ can be
calculated from the long exact sequence of the pair
$(M^2,M^2\setminus \Delta_M)$,
\begin{equation}\label{eqn:long-pair}
\ldots \longrightarrow H^\ast(M^2,M^2\setminus \Delta_M)
\stackrel{\alpha}{\longrightarrow} H^\ast(M^2)
\stackrel{\beta}{\longrightarrow} H^\ast(F_2(M))
\longrightarrow\ldots
\end{equation}
We are interested in the (dual) Stiefel-Whitney classes
(Proposition~\ref{prop:dual}) so we tacitly assume that all
cohomology has coefficients $\mathbb{F}_2$ unless otherwise noted.

By naturality, in order to check non-triviality of
$\overline{w}_k(T(F_2(M)))$, it is sufficient to check that
$\overline{w}_k(M^2)$ is not in the image of the map $\alpha$.

\medskip
The image $A:={\rm Image}(\alpha)$ of $\alpha$ is determined in
\cite[Theorem~11.11]{Mi-Sta} (see also \cite[Chapter VI,
Section~12]{Bredon}). It is generated, as a $H^\ast(M)$-module, by
the ``diagonal cohomology class''
\begin{equation}\label{eqn:generator}
u'' = \sum_{i=1}^r b_i\times b_i^\sharp
\end{equation}
where $\{b_i\}_{i=1}^r$ is an additive basis of $H^\ast(M)$ and
$b_i^\sharp$ the class dual to $b_i$.

\medskip
There are two actions of the ring $H^\ast(M)$ on $H^\ast(M\times
M)$ associated with the projections $\pi_1,\pi_2 : M^2\rightarrow
M$. However in light of \cite[Lemma~11.8]{Mi-Sta}, which says that
if $a\in H^\ast(M)$ then
$$
(a\times 1)\cup u'' = (1\times a)\cup u'',
$$
these two actions have the same effect on $u''$. As a consequence
we obtain the following proposition.

\begin{prop}\label{prop:image}
\begin{equation}
\begin{array}{ccccccc}
A & = & \mbox{\rm Image}(\alpha) & = & H^\ast(M)\cdot u''& = &
\{(1\times a)\cup u'' \mid a\in H^\ast(M)\} \\
&&&&& = & \{(a\times 1)\cup u'' \mid a\in H^\ast(M)\}
\end{array}
\end{equation}
\end{prop}

\medskip
The following proposition provides a simple and efficient
criterion for testing if a class is in the image of the map
$\alpha$. Note that the condition $k\leq n-1$ is essential since
\begin{equation}
H^{2n}(M\times M) \cong H^n(M)\otimes H^n(M)\subset \mbox{\rm
Image}(\alpha).
\end{equation}
\begin{prop}\label{prop:sinisa}
\label{diagonal} Let $M$ be a closed and smooth $n$-dimensional
manifold. Let $k\leq n-1$ and assume that $\theta\in H^k(M)\otimes
H^k(M)\subset H^\ast(M\times M)\cong H^\ast(M)\otimes H^\ast(M)$
is a non-zero class. Then $\theta\notin {\rm Image}(\alpha )$.
More generally, if $\omega\in H^{n+p}(M\times M)$ is a non-zero
class which is in the image of $\alpha$ then it must have as a
component of bidegree $(p,n)$ a non-zero ``edge class'' of the
form $a\times z$ for some $a\in H^p(M)$, where $z\in H^n(M)$ is
the fundamental cohomology class of $M$.
\end{prop}
\medskip

\noindent {\bf Proof:} If $z\in H^n(M)$ is the fundamental
cohomology class of $M$ then the diagonal class $u''$ has the
following form
$$u''=z\times 1 + x_1\times y_{1}+\cdots +x_r\times y_r + 1\times
z$$ where $x_i\times y_i$ is a class of bidegree $(t,n-t)$ for
some $0<t<n$. If $\omega\in H^{n+p}(M\times M)$ is in the image of
$\alpha$ then we deduce from Proposition~\ref{prop:image} that
$$\omega = (a\times 1)u'' = az\times 1 + A + a\times z$$
where $A=ax_1\times y_1 +\cdots$ is a class whose homogeneous
components are of bidegree $(q,n+p-q)$ for some $q>p$, and the
proposition follows. \hfill $\square$

\begin{cor}\label{cor:lepo-1}
If $k:={\rm max}\{i\mid \overline{w}_i(M)\neq 0\}$ then
$\overline{w}_{2k}(T(F_2(M)))=w_{2k}(\nu)\neq 0$.
\end{cor}

\medskip\noindent
{\bf Proof:} It follows from the naturality of Stiefel-Whitney
classes that
$$
w_{2k}(\nu) = \overline{w}_{2k}(T(F_2(M))) =
\beta(\overline{w}_{2k}(M^2)) = \beta(\overline{w}_k(M)\times
\overline{w}_k(M)).
$$
We observe that $k\leq n-1$. Indeed, each $n$-dimensional smooth
manifold can be embedded in $\mathbb{R}^{2n}$ and
$\overline{w}_n(M)=0$ by \cite[Corollary~11.4.]{Mi-Sta}.

Since $k\leq n-1$ we are allowed to use
Proposition~\ref{prop:sinisa} which implies that
$\overline{w}_k(M)\times \overline{w}_k(M)\notin {\rm
Image}(\alpha)$. From here and the exactness of the sequence
(\ref{eqn:long-pair}) we finally deduce that $w_{2k}(\nu)\neq 0$.

\hfill $\square$

\section{Real projective spaces}

As a first application let us analyze the case when $M =
\mathbb{R}P^n$ is the $n$-dimensional real projective space.

\medskip
The cohomology algebra $$H^\ast(\mathbb{R}P^n)\cong
\mathbb{F}_2[t]/(t^{n+1}=0)$$ is a truncated polynomial ring with
one generator $t\in H^1(\mathbb{R}P^n)$.

\medskip
The total Stiefel-Whitney class of $T(\mathbb{R}P^n)$ is given
(\cite[Theorem~4.5.]{Mi-Sta}) by the formula
\begin{equation}\label{eqn:total-class-proj}
w(\mathbb{R}P^n) = (1+t)^{n+1}
\end{equation}
and the dual classes are
$$
\overline{w}(\mathbb{R}P^n) = w(\mathbb{R}P^n)^{-1}.
$$

Suppose that $n=2^r$ is a power of $2$. Then
$$
w(\mathbb{R}P^n)= 1 +t + t^n \quad\mbox{\rm and}\quad
\overline{w}(\mathbb{R}P^n) = 1+t+t^2+\cdots + t^{n-1}.
$$
It follows that

$$\overline{w}_{n-1}(\mathbb{R}P^n)=t^{n-1}\neq 0.$$

As a consequence of Corollary~\ref{cor:lepo-1} we obtain that
$\overline{w}_{2n-2}(F_2(\mathbb{R}P^n))\neq 0$ and deduce from
Proposition~\ref{prop:dual} the following result.

\begin{prop}\label{prop:proj-donja-ocena}
If $f : \mathbb{R}P^n\rightarrow \mathbb{R}^N$ is a totally skew
embedding and $n=2^r$ for some $r$ then
$$
N \geq 4n-1.
$$
\end{prop}

\begin{cor}\label{cor:proj} For each integer $n$,
$$
N(\mathbb{R}P^n)\geq 4\cdot 2^{[\log_2 n]}-1.
$$
\end{cor}

It follows from Proposition~\ref{prop:proj-donja-ocena} and the
inequalities (\ref{eqn:lower-upper}) that if $n=2^r$ is a power of
$2$ then $N(\mathbb{R}P^n)$ is $4n-1, 4n$ or $4n+1$, in particular
$N(\mathbb{R}P^2)$ is $7,8$ or $9$.

\section{Products of real projective spaces}
\label{sec:prod-proj-spaces}

Suppose that $X=\mathbb{R}P^{n_1}\times \cdots \times
\mathbb{R}P^{n_k}$ is a product of real projective spaces where
each $n_i=2^{r_i}$ is a power of $2$. Let $n={\rm dim}(X) =
n_1+\cdots +n_k$.

\medskip The cohomology $H^*(X)\cong
\mathbb{F}_2[u_1,...,u_k]/(u_1^{n_1+1}=\ldots =u_k^{n_k+1}=0)$ of
$X$ is a truncated polynomial algebra with $k$ generators
$u_1,...,u_k\in H^1(X)$. Since $T(X)=T(\mathbb{R}P^{n_1})\times
\cdots \times T(\mathbb{R}P^{n_k}),$ the total Stiefel-Whitney
class of $T(X)$ is given by the formula
$$w(X)=(1+u_1)^{n_1+1}\cdots (1+u_k)^{n_k+1},$$

\noindent and its dual total class is $\overline{w}(X)=w(X)^{-1}$.

By assumption all integers $n_i$ are powers of $2$, hence
$$w(X)=(1+u_1+u_1^{n_1})\cdots (1+u_k+u_k^{n_k}),$$

\noindent and the dual class has the form
$$\overline{w}(X)=(1+u_1+\cdots +u_1^{n_1-1})
\cdots (1+u_k+\cdots +u_k^{n_k-1}).$$

From here we deduce that $\overline{w}_{n-k}=u_1^{n_1-1}\cdots
u_k^{n_k-1}$ is non-zero and observe, by a reference to
Proposition~\ref{prop:sinisa} and Corollary~\ref{cor:lepo-1}, that
$\overline{w}_{2n-2k}(F_2(X)) \neq 0$. This fact allows us to use
Proposition~\ref{prop:dual} which in turn implies the following
theorem.

\begin{theo}\label{thm:proj-product} Suppose that
$X=\mathbb{R}P^{n_1}\times \cdots \times \mathbb{R}P^{n_k}$ where
$n_i=2^{r_i}$ are powers of $2$. Let $n={\rm dim}(X)=n_1+\cdots
+n_k$. If there exists a totally skew embedding of $X$ in
$\mathbb{R}^N$ then $N\geq 4n-2k+1$. In particular if $n_i\neq
n_j$ for $i\neq j$,
$$
N(X)\geq 4n-2\alpha (n)+1
$$
where $\alpha(n)$ is the number of non-zero digits in the binary
representation of $n$.
\end{theo}

\section{Complex manifolds}

In some cases, for example if  $M$ is a complex manifold, it may
be convenient to use Pontryagin classes for estimating the
invariant $N(M)$. However, the inequalities obtained by the use of
Pontryagin classes are in general not as sharp as the inequalities
obtained with the aid of Stiefel-Whitney classes so we focus on
the latter method.

\subsection{Complex projective spaces}
\label{sec:compl-proj}

The cohomology of the complex projective space with $\mathbb{Z}$
coefficients is a truncated polynomial algebra so by the Universal
Coefficient Theorem we have $H^\ast(\mathbb{C}P^n;
\mathbb{F}_2)\cong \mathbb{F}_2[t]/(t^{n+1}=0)$ where $\mbox{\rm
deg}(t)=2$.

Since the second Stiefel-Whitney class $w_2$ of any oriented
$2$-plane bundle is the mod 2 reduction of the Euler class, we
observe that $t = w_2(\xi_{\mathbb{R}}) =
w_2(\xi^\ast_{\mathbb{R}})$ where $\xi$ is the canonical complex
line bundle over $\mathbb{C}{P}^n$ and $\xi_{\mathbb{R}}$ the
underlying real 2-plane bundle.

\medskip
The complex tangent bundle of the projective space $\mathbb{C}P^n$
is
$$
T(\mathbb{C}P^n)\cong {\rm Hom}(\xi,\xi^\perp)
$$
where $\xi^\perp$ is the complex $n$-plane bundle, complementary
to the tautological complex line bundle $\xi$. Since ${\rm
Hom}(\xi,\xi)\cong \varepsilon^1_{\mathbb{C}}$ is a trivial
complex line bundle, we conclude that
\begin{equation}\label{eqn:complex-proj-tang-bundle}
T(\mathbb{C}P^n)\oplus\varepsilon^1_{\mathbb{C}} \cong
(\xi^\ast)^{\oplus (n+1)}
\end{equation}
where $\xi^\ast$ is the line bundle dual to $\xi$. By forgetting
the complex structure (realification) we obtain the isomorphism of
real bundles
\begin{equation}\label{eqn:complex-proj-real-bundle}
T(\mathbb{C}P^n)_{\mathbb{R}}\oplus\varepsilon^2_{\mathbb{R}}
\cong (\xi^\ast_{\mathbb{R}})^{\oplus (n+1)}.
\end{equation}
It follows that the total Stiefel-Whitney class of
$T(\mathbb{C}P^n)_{\mathbb{R}}$ is
$$
w = w(T(\mathbb{C}P^n)_{\mathbb{R}}) = (1+t)^{n+1} = (1+w_2)^{n+1}
$$
where $w_2 = w_2(\xi_{\mathbb{R}})\cong
w_2(\xi^\ast_{\mathbb{R}})=t$ is the second Stiefel-Whitney class
of the realification of the canonical bundle $\xi$.

\medskip Consequently, the dual Stiefel-Whitney class is
\begin{equation}\label{eqn:dual-complex}
\overline{w} = (1+ w_2)^{-n-1}=\sum_{j=0}^n {{n+j}\choose{j}}
w_2^j.
\end{equation}

We observe that the top class $\overline{w}_{2n}$ is always zero,
which is an instance of a much more general result of Massey
(Theorem~\ref{thm:Massey}). Following Corollary~\ref{cor:lepo-1}
we search for the largest value of $j$ such that
$\overline{w}_{2j}={{n+j}\choose{j}}w_2^j\neq 0$. We observe that
$\overline{w}_{2n-2}\neq 0$ precisely when $n=2^r$ is a power of
$2$.

Again, by invoking Corollary~\ref{cor:lepo-1},  we conclude that
$\overline{w}_{4n-4}(\nu)\neq 0$ and finally by
Proposition~\ref{prop:dual} we obtain the following result.

\begin{theo}\label{thm:compl-proj}
Suppose that $n=2^r$ for some $r\geq 0$. Then
$$
N(\mathbb{C}P^n) \geq 4n + (4n-4) + 1 = 4\cdot\mbox{\rm
dim}(\mathbb{C}P^n) - 3.
$$
\end{theo}

\subsection{Products of complex projective spaces}
\label{sec:prod-compl-proj}

Suppose that $X = \mathbb{C}P^{n_1}\times\cdots\times
\mathbb{C}P^{n_k}$. As in Section~\ref{sec:prod-proj-spaces} we
focus on the case when $n_i=2^{r_i}$ for some $r_i$. As before $n
= (1/2){\rm dim}(X)=n_1+\cdots +n_k$. The cohomology ring of the
space $X$ with $\mathbb{F}_2$ coefficients is
$$H^*(X)\cong \mathbb{F}_2[u_1,...,u_k]/
(u_1^{n_1+1}=\ldots =u_k^{n_k+1}=0)$$ where ${\rm deg}(u_1)=\ldots
= {\rm deg}(u_{n_k})=2$.

\medskip
We have already observed in Section~\ref{sec:compl-proj} that if
$n=2^r$ then
\begin{equation}\label{eqn:dual-complex-jos}
\overline{w}(T(\mathbb{C}P^n)) = \sum_{j=0}^n {{n+j}\choose{j}}
t^j = 1 + t + \cdots + t^{n-1}.
\end{equation}

It follows, as in Section~\ref{sec:prod-proj-spaces}, that the
total dual Stiefel-Whitney class of $T(X)$ has the form
$$\overline{w}(X)=(1+u_1+\cdots +u_1^{n_1-1})\cdots
(1+u_k+\cdots +u_k^{n_k-1}).$$

We conclude that $\overline{w}_{2n - 2k}= u_1^{n_1-1}\cdots
u_k^{n_k-1}$ is non-trivial, and as a consequence of
Proposition~\ref{prop:dual} and Corollary~\ref{cor:lepo-1} obtain
the following estimate.

\begin{theo}\label{thm:compl-proj-product}
Suppose that $X=\mathbb{C}P^{n_1}\times \cdots \times
\mathbb{C}P^{n_k}$ where $n_i=2^{r_i}$ are powers of $2$ and let
$n={\rm dim}_{\mathbb{C}}(X) = (1/2){\rm
dim}_{\mathbb{R}}(X)=n_1+\cdots +n_k$. Then $N(X)\geq 8n-4k+1$. In
particular if all integers $n_i$ are distinct,
$$
N(X)\geq 8n-4\alpha (n)+1 = 4\cdot{\rm dim}_{\mathbb{R}}(X) -
4\alpha(n) +1
$$
where $\alpha(n)$ is the number of non-zero digits in the binary
representation of $n$.
\end{theo}

\section{Grassmannians}\label{sec:grassmannians}

We illustrate our method also for some cases of the Grassmann
manifold $G_k(\mathbb{R}^{n+k})$ of $k$-dimensional subspaces of
$\mathbb{R}^{n+k}$, and some cases of the oriented Grassmann
manifold $\tilde{G}_k(\mathbb{R}^{n+k})$ of oriented
$k$-dimensional subspaces of $\mathbb{R}^{n+k}$.

\medskip
Let $\gamma_k$ be the canonical vector bundle over
$X=G_k(\mathbb{R}^{n+k})$, and $\tau$ the tangent bundle. Then
from the relation $\tau \oplus {\rm Hom}(\gamma_k,\gamma_k)\simeq
(n+k)\gamma_k$ we obtain
\begin{equation}\label{eqn:1}
w(X)\cdot w(\gamma_k\otimes \gamma_k^*)=w(\gamma_k)^{n+k},
\end{equation}
where $w(\gamma_k\otimes \gamma_k^*)=p_k(w_1,...,w_k)$, $p_k$ is
the polynomial over $\mathbb{F}_2$ defined by
\begin{equation}\label{eqn:2}
p_k(\sigma_1,...,\sigma_k)=\prod_{i=1}^k\prod_{j=1}^k(1+x_i+x_j),
\end{equation}
and $\sigma_1,...,\sigma_k$ are the elementary symmetric
polynomials in variables $x_1,...,x_k$, see \cite[Problem
7C]{Mi-Sta}.

In the special case when $k=2$ and $k=3$, by a direct computation
we check that $w(\gamma_2 \otimes \gamma_2^*)=1+w_1^2$ and
$w(\gamma_3\otimes \gamma_3^*)= 1+w_1^4+w_2^2+w_1^2w_2^2+w_3^2$.

Since $w(\gamma_k)=1+w_1+w_2+\cdots +w_k$, it follows that the
total Stiefel-Whitney class of the complementary bundle to the
tangent bundle $\tau$ equals
\begin{equation}\label{eqn:3}
\overline{w}(X)=w(\gamma_k\otimes \gamma_k^*)(1+w_1+w_2+\cdots
+w_k)^{-(n+k)}.
\end{equation}

Completely analogous formulae are true in the case of the oriented
Grassmann manifold $\tilde{X}=\tilde{G}_k(\mathbb{R}^{n+k})$, the
only difference being the vanishing of the first Stiefel-Whitney
class, $w_1(\tilde{\gamma}_k)=0$.

\subsection{$G_k(\mathbb{R}^{n+k})$}

First we treat the case $k=2$ and $n=2^r$, that is the case of the
Grassmann manifold $G_2\left(\mathbb{R}^{2^r+2}\right)$. We shall
need the following lemma.

\begin{lema}
\label{l1} The class $w_1^2w_2^{2^r-2}\in
H^{2^{r+1}-2}\left(G_2\left(\mathbb{R}^{2^r+2}\right)\right)$ is
non-trivial.
\end{lema}
\medskip

\noindent {\bf Proof:} Let us assume, to the contrary, that
$w_1^2w_2^{2^r-2}=0$. Then $w_1^2w_2^{2^r-1}=0$. Since the map
\begin{equation}\label{eqn:4}
H^{2^{r+1}-1}\left(G_2\left(\mathbb{R}^{2^r+2}\right)\right)
\stackrel{\cup w_1} {\longrightarrow}
H^{2^{r+1}}\left(G_2\left(\mathbb{R}^{2^r+2}\right)\right)
\end{equation}
is an isomorphism by Poincar\' e duality, it follows that
$w_1w_2^{2^r-1}=0$.

Let us show that $w_1^{2^{r+1}-2}$ and $w_2^{2^r-1}$ are
non-trivial classes in
$H^{2^{r+1}-2}\left(G_2\left(\mathbb{R}^{2^r+2}\right)\right)$.
The first observation is a consequence of a result of Stong
\cite{s} about the height of $w_1$, which is in this case
$\mbox{{\rm ht}}(w_1)=2^{r+1}-2$. The second observation follows
from the well-known fact that $w_k^n$ is a non-trivial element in
$H^{kn}(G_k(\mathbb{R}^{k+n}))$. Let us show that these two
classes are different. We have
\begin{equation}\label{eqn:5}
Sq^2\left(w_1^{2^{r+1}-2}\right)={2^{r+1}-2 \choose
2}w_1^{2^{r+1}}=0,
\end{equation}
again by the same result of Stong. Since by the Wu formula
$Sq^1(w_2)=w_1w_2$ (see \cite[Problem 8A]{Mi-Sta}), we have
\begin{equation}\label{eqn:6}
Sq^2\left(w_2^{2^r-1}\right)=(2^r-1)w_2^{2^r}+{2^r-1 \choose 2}
w_1^2w_2^{2^r-1}=w_2^{2^r}\neq 0.
\end{equation}
So, $H^{2^{r+1}-2}\left(G_2\left(\mathbb{R}^{2^r+2}\right)\right)
\cong \mathbb{Z}/2\oplus \mathbb{Z}/2$ is generated by
$w_1^{2^{r+1}-2}$ and $w_2^{2^r-1}$.

The map $\phi :
H^{2^{r+1}-2}\left(G_2\left(\mathbb{R}^{2^r+2}\right)\right)
\stackrel{\cup w_1} {\longrightarrow}
H^{2^{r+1}-1}\left(G_2\left(\mathbb{R}^{2^r+2}\right)\right)$
satisfies the relations $\phi (w_1^{2^{r+1}-2})=w_1^{2^{r+1}-1}=0$
and $\phi (w_2^{2^r-1})=w_1w_2^{2^r-1}=0$, as we proved in the
beginning of the proof. So, $\phi =0$. This is a contradiction,
since $2^{r+1}-1$ is odd and
$H^{2^{r+1}-1}\left(G_2\left(\mathbb{R}^{2^r+2}\right)\right)
\cong \mathbb{Z}/2$ could be generated only by the element of the
type $w_1^{2s-1}w_2^t$. So, our assumption is false, and we have
$w_1^2w_2^{2^r-2}\neq 0$. \hfill $\square$
\bigskip

\begin{theo}
\label{t1} $N\left(G_2\left(\mathbb{R}^{2^r+2}\right)\right)\geq
4\cdot 2^{r+1}-3.$
\end{theo}
\medskip

\noindent {\bf Proof:} The total Stiefel-Whitney class of the
normal bundle of $X=G_2\left(\mathbb{R}^{2^r+2}\right)$ equals, by
the equation (\ref{eqn:3}),
\begin{eqnarray*}
\overline{w}(X) & = & (1+w_1^2)(1+w_1+w_2)^{-(2^r+2)}\\
& = & (1+w_1^2)(1+w_1+w_2)^{-2^{r+1}}(1+w_1+w_2)^{2^r-2}\\
& = & (1+w_1^2)(1+w_1^{2^{r+1}})(1+w_1+w_2)^{2^r-2}\\
& = & (1+w_1^2)(1+w_1+w_2)^{2^r-2}.
\end{eqnarray*}

It follows that $\overline{w}_{2^{r+1}-2}=w_1^2w_2^{2^r-2}\neq 0$,
which by Corollary~\ref{cor:lepo-1} and
Proposition~\ref{prop:dual} implies the inequality
$$N\left(G_2\left(\mathbb{R}^{2^r+2}\right)\right)\geq 4\cdot 2^{r+1}-3=4\cdot \dim
G_2\left(\mathbb{R}^{2^r+2}\right)-3.$$ \hfill $\square$
\bigskip

As an illustration of our methods in the case $k>2$, we outline
the computations in the particular case of the Grassmann manifold
$G_3(\mathbb{R}^7)$.

\begin{theo}
\label{t2} $N(G_3(\mathbb{R}^7))\geq 43.$
\end{theo}
\medskip

\noindent {\bf Proof:} The cohomology of $X=G_3(\mathbb {R}^7)$ is
generated by the Stiefel-Whitney classes $w_1,w_2,w_3$ subject to
the relation
$(1+w_1+w_2+w_3)(1+\overline{w}_1+\overline{w}_2+\overline{w}_3
+\overline{w}_4)=1$. It follows that $\overline{w}_1=w_1$,
$\overline{w}_2=w_1^2+w_2$, $\overline{w}_3=w_1^3 +w_3$, and
$\overline{w}_4=w_1^4+w_1^2w_2+w_2^2$. Moreover, $R_1:=
w_1^5+w_1w_2^2+w_1^2w_3=0$, $R_2:= w_1^4w_2+w_1^2w_2
^2+w_2^3+w_1^3w_3+w_3^2=0$, and $R_3:=
w_1^4w_3+w_1^2w_2w_3+w_2^2w_3=0$.

It could be shown that a consequence of these relations is also
the relation
$$0=(w_1^3+w_1w_2+w_3)R_1+w_1^2R_2+w_1R_3=w_1^8.$$

It requires a few more steps to show that the class
$w_1^2w_2^2w_3+w_3^3=w_1^5w_2^2$ is non-trivial. It is actually
one of the additive generators of the cohomology group
$H^9(G_3(\mathbb{R}^7))$.

The total Stiefel-Whitney class of the bundle complementary to the
tangent bundle equals
\begin{eqnarray*}
\overline{w}(X) & = &
(1+w_1+w_2+w_3)^{-7}\cdot (1+w_1^4+w_2^2+w_1^2w_2^2+w_3^2)\\
 & = & (1+w_1+w_2+w_3)(1+w_1+w_2+w_3)^{-8}(1+w_1^4+w_2^2+w_1^2w_2^2+w_3^2)\\
 & = & (1+w_1+w_2+w_3)(1+w_1^4+w_2^2+w_1^2w_2^2+w_3^2).
\end{eqnarray*}

We already noticed that $\overline{w}_9(X)=w_1^2w_2^2w_3+w_3^3$ is
a non-trivial class, and it is the top-dimensional one.

Altogether, we conclude that $N(G_3(\mathbb{R}^7)) \geq
24+1+18=43.$ \hfill $\square$
\bigskip

Let us add that in a similar way but more easily one obtains by
the same method also: $N(G_2(\mathbb{R}^5))\geq 21$,
$N(G_2(\mathbb{R}^7))\geq 29$, $N(G_3(\mathbb{R}^6))\geq 31$, and
$N(G_3(\mathbb{R}^8))\geq 43$.

\subsection{$\tilde{G}_k(\mathbb{R}^{n+k})$}

For comparison we include an analysis of some cases where the
manifold $M$ is the Grassmannian of all oriented $k$-dimensional
subspaces in $\mathbb{R}^{n+k}$.

\medskip Let us denote by $p :
\tilde{G}_k(\mathbb{R}^{n+k})\rightarrow G_k(\mathbb{R}^{n+k})$
the two-fold covering. Then, $\tilde{w}_i=
w_i(\tilde{G}_k(\mathbb{R}^{n+k}))=p^*(w_i(G_k(\mathbb{R}^{n+k})))$,
and we know that $\tilde{w}_1=0$. Since $\tilde{w}_i=p^*(w_i)\neq
0$ implies $w_i\neq 0$, the estimates obtained in this way for the
oriented Grassmann manifold $\tilde{G}_k(\mathbb{R}^{n+k})$ cannot
be better than those for $G_k(\mathbb{R}^{n+k})$.

However, the cohomology ring of the oriented Grassmann manifold
$H^*(\tilde{G}_k(\mathbb{R}^{n+k}))$ is more complicated, and it
is more difficult to determine which Stiefel-Whitney classes are
non-trivial in this case. Aside from triviality of $\tilde{w}_1$,
we know that $H^*(\tilde{G}_k(\mathbb{R}^{n+k}))$ has some
additional generators and some additional relations.

\medskip
Let $B_n^k = (w_1)$ be the principal ideal in
$H^\ast(G_2(\mathbb{R}^{n+k}))$ generated by $w_1$. In order to
determine which Stiefel-Whitney classes are non-trivial in the
oriented case, we use the calculations in
$H^*(G_k(\mathbb{R}^{n+k}))$ and the Gysin exact sequence in
cohomology (cf.\  \cite[Theorem 12.4]{Mi-Sta}),
$$\cdots \rightarrow H^{i-1}(G_k(\mathbb{R}^{n+k}))\stackrel{\cup w_1}
{\longrightarrow} H^{i}(G_k(\mathbb{R}^{n+k})) \stackrel{p^*}
{\longrightarrow} H^{i}(\tilde{G}_k(\mathbb{R}^{n+k})) \rightarrow
\cdots .
$$

From the exactness of this sequence we deduce that for a given
Stiefel-Whitney class $w_{i_1}^{j_1}\cdots w_{i_r}^{j_r}\in
H^*(G_k(\mathbb{R}^{n+k}))$,  $p^*(w_{i_1}^{j_1}\cdots
w_{i_r}^{j_r})=0$ if and only if $w_{i_1}^{j_1}\cdots
w_{i_r}^{j_r}\in B_n^k = (w_1)$.

Also, we easily check that in this case the polynomials $p_2$ and
$p_3$ from the equation (\ref{eqn:2}) reduce to the following,
$\tilde{p}_2(\tilde{w}_2)=1$ and
$\tilde{p}_3(\tilde{w}_2,\tilde{w}_3)
=1+\tilde{w}_2^2+\tilde{w}_3^2$.

\medskip
We now turn to the case $k=2$. Let us determine the height of the
class $\tilde{w}_2$ in $H^*(\tilde{G}_2(\mathbb{R}^{n+2}))$,
$\textrm{ht}(\tilde{w}_2)=\max \{ m\in \mathbb{N} \mid
\tilde{w}_2^m\neq 0\}$.

In $H^{*}(G_2(\mathbb{R}^{n+2}))$ we have $(1+w_1+w_2)(1+
\overline{w}_1+...+\overline{w}_n)=1$, and so
\begin{equation}\label{eqn:jos-malo}
\overline{w}_r=w_1\overline{w}_{r-1}+w_2\overline{w}_{r-2}, \quad
3\leq r\leq n.
\end{equation}

If as before $B_n^2=(w_1)$ is the principal ideal in
$H^{*}(G_2(\mathbb{R}^{n+2}))$ generated by $w_1$, then
inductively, using the relations (\ref{eqn:jos-malo}),  we show
that
$$\overline{w}_{2k-1}\in B_n^2, \quad 2k-1\leq n$$ and
$$\overline{w}_{2k}\equiv w_2^k \quad (\textrm{mod }B_n^2), \quad 2k\leq n.$$
Note that in dimensions $\leq n$ there are no polynomial relations
among $w_1$ and $w_2$.

\begin{lema}
\label{l2} $\mbox{{\rm ht}}(\tilde{w}_{2})=[\frac{n}{2}]$.
\end{lema}
\medskip

\noindent {\bf Proof:} It is well known that $\ker p^{*}=B_{n}^2$.
The dimension of $w_2^{[\frac{n}{2}]}$ is $2\cdot
[\frac{n}{2}]\leq 2\cdot \frac{n}{2}=n$, hence this class cannot
be written as a multiple of $w_1$ (for in dimensions $\leq n$
there are no relations among $w_1$ and $w_2$). Thus
$w_2^{[\frac{n}{2}]} = \notin \ker p^*$ and so
$\tilde{w}_2^{[\frac{n}{2}]} = p^*(w_2^{[\frac{n}{2}]}) \neq 0$.
\smallskip

In order to show that $\tilde{w}_2^{[\frac{n}{2}]+1}=0$ we
distinguish two cases.

If $n=2l$, in dimension $2l+2$ we have the relation
$w_2\overline{w}_{2l}=0$. But $\overline{w}_{2l}=w_2^l+w_1\cdot u$
for some class $u$, so $w_2^{l+1}+w_1w_2u=0$ and $w_2^{l+1}\in
B_n^2$. So we obtain
$$\tilde{w}_2^{[\frac{n}{2}]+1}=\tilde{w}_2^{l+1}=p^*(w_2^{l+1})=0.$$

If $n=2l+1$, in dimension $2l+2$ we have the relation
$w_1\overline{w}_ {2l+1}+w_2\overline{w}_{2l}=0$. The first
summand belongs to $B_{n}^2$, so we show (as in the first case)
that $w_2^{l+1}\in B_n^2$. Here we also have that
$l=[\frac{n}{2}]$ and the Lemma follows. \hfill $\square$
\medskip

Let us also notice that by the equation (\ref{eqn:3}) and the fact
that $\tilde{p}_2=1$, the total Stiefel-Whitney class of the
complementary normal bundle to the tangent bundle of the space
$X=\tilde{G}_2(\mathbb{R}^{n+2})$ equals
$$\overline{w}(X)=(1+\tilde{w}_2)^{-(n+2)}=((1+\tilde{w}_2)^{-1})^{n+2}.$$

In the light of Lemma \ref{l2}, we have
$$(1+\tilde{w}_2)\left(1+\tilde{w}_2+\tilde{w}_2^2+\cdots +\tilde{w}_2^{[\frac{n}{2}]}\right)=
1+\tilde{w}_2^{[\frac{n}{2}]+1}=1,$$

\noindent and so
$(1+\tilde{w}_2)^{-1}=1+\tilde{w}_2+\tilde{w}_2^2+\cdots
+\tilde{w}_2^{[\frac{n}{2}]}$.

Finally, we obtain
$$\overline{w}(\tilde{G}_2(\mathbb{R}^{n+2}))=
\left(1+\tilde{w}_2+\tilde{w}_2^2+\cdots
+\tilde{w}_2^{[\frac{n}{2}]}\right)^{n+2}.$$

\begin{theo}
\label{t3} $N(\tilde{G}_2(\mathbb{R}^{2^r+2}))\geq 3\cdot
2^{r+1}+1= 3\cdot \dim \tilde{G}_2(\mathbb{R}^{2^r+2})+1.$
\end{theo}
\medskip

\noindent {\bf Proof:} Substituting $n=2^r$ in the above
considerations and using Lemma \ref{l2}, we have
\begin{eqnarray*}
\overline{w}(\tilde{G}_2(\mathbb{R}^{2^r+2})) & = &
(1+\tilde{w}_2+\tilde{w}_2^2+\cdots
+\tilde{w}_2^{2^{r-1}})^{2^r+2}\\
& = & (1+\tilde{w}_2+\tilde{w}_2^2+\cdots
+\tilde{w}_2^{2^{r-1}})^2\\
& = & 1+\tilde{w}_2^2+\tilde{w}_2^4+\cdots +\tilde{w}_2^{2^{r-1}}.
\end{eqnarray*}

By Lemma \ref{l2}, $\tilde{w}_2^{2^{r-1}}\neq 0$, and by
Corollary~\ref{cor:lepo-1},
$N(\tilde{G}_2(\mathbb{R}^{2^r+2}))\geq 1+2\cdot 2^{r+1}+2\cdot
2^r=3\cdot 2^{r+1}+1=3\cdot \dim
\tilde{G}_2(\mathbb{R}^{2^r+2})+1. $ \hfill $\square$

\medskip
Let us add that by the same methods one easily obtains
$N(\tilde{G}_2(\mathbb{R}^{2^r+1}))\geq 3\cdot 2^{r+1}-7=3\cdot
\dim \tilde{G}_2(\mathbb{R}^{2^r+1})-1$,
$N(\tilde{G}_2(\mathbb{R}^{2^r+3}))\geq 3\cdot 2^{r+1}+5=3\cdot
\dim \tilde{G}_2(\mathbb{R}^{2^r+3})-1$, and
$N(\tilde{G}_2(\mathbb{R}^{2^r+4}))\geq 3\cdot 2^{r+1}+9=3\cdot
\dim \tilde{G}_2(\mathbb{R}^{2^r+4})-3$. It is also seen from the
proof that our method cannot give better lower bounds in all these
cases.

Let us now prove the result corresponding to Theorem \ref{t2} in
the case of the oriented Grassmannian.

\begin{theo}
\label{t4} $N(\tilde{G}_3(\mathbb{R}^7))\geq 41.$
\end{theo}
\medskip

\noindent {\bf Proof:} In the cohomology of the oriented
Grassmannian $\tilde{X}=\tilde{G}_3(\mathbb{R}^7)$ we have
\begin{eqnarray*}
\overline{w}(\tilde{X})=p^*(\overline{w}(X)) & = &
(1+\tilde{w}_2+\tilde{w}_3)^{-7}(1+\tilde{w}_2^2+\tilde{w}_3^2)\\
& = & (1+\tilde{w}_2+\tilde{w}_3)(1+\tilde{w}_2+\tilde{w}_3)^{-8}
(1+\tilde{w}_2^2+\tilde{w}_3^2)\\
& = & (1+\tilde{w}_2+\tilde{w}_3)(1+\tilde{w}_2^2+\tilde{w}_3^2).
\end{eqnarray*}

In the previous subsection we showed that the class
$\overline{w}_9(X)=w_1^2w_2^2w_3+w_3^3=w_1^5w_2^2$ is non-trivial,
but it is trivial in the cohomology of the oriented Grassmannian.
However, it can be shown (using the computations in $H^*(G_3
(\mathbb{R}^7))$) that the class
$\overline{w}_8(X)=w_1^2w_2^3+w_2w_3^2$ cannot be written as a
product of $w_1$ with some other class. So,
$$
\overline{w}_8(\tilde{X})=\tilde{w}_2\tilde{w}_3^2\neq 0.
$$

As a consequence, by Corollary~\ref{cor:lepo-1} we have
$N(\tilde{G}_3(\mathbb{R}^7))\geq 1+24+16=41. \hfill\square$

\bigskip
We end this section by presenting a slightly different method of
calculation applied to the Grassmannian
$Y=\tilde{G}_3(\mathbb{R}^{13})$. In this case $\dim (Y)=30$ and
\begin{eqnarray*}
\overline{w}(Y) & = &
(1+\tilde{w}_2^2+\tilde{w}_3^2)(1+\tilde{w}_2+\tilde{w}_3)^{-13}\\
& = & (1+\tilde{w}_2^2+\tilde{w}_3^2)(1+\tilde{w}_2+\tilde{w}_3)^3
(1+\tilde{w}_2+\tilde{w}_3)^{-16}\\
& = &
(1+\tilde{w}_2^2+\tilde{w}_3^2)(1+\tilde{w}_2+\tilde{w}_3)^3\\
& = & \tilde{w}_3^5+\tilde{w}_2\tilde{w}_3^4+\cdots ,
\end{eqnarray*}

\noindent where dots represent some  lower-dimensional classes. In
order to check whether some of the classes
$\overline{w}_{15}(Y)=\tilde{w}_3^5$ and
$\overline{w}_{14}(Y)=\tilde{w}_2\tilde{w}_3^4$ are non-trivial we
use a criterion from \cite{Kor06}. It says that
$w_2^{i_2}w_3^{i_3}\in H^*(G_3(\mathbb{R}^{13}))$ cannot be
expressed as a multiple of the class $w_1$ if it does not belong
to the ideal $J_{n,3}$ of $\mathbb{Z}_2[w_2,w_3]$ generated by the
homogeneous components of
$$
\frac 1{1+w_2+w_3}=(1+w_2+w_3)^{2^{r+3}-1}=
\sum_{i=0}^{2^{r+3}-1}\sum_{j=0}^i{i \choose j}w_2^{i-j}w_3^j
$$

\noindent in dimensions $n-2,n-1,n$, which we respectively denote
by $g_{n-2},g_{n-1},g_n$. The integer $r$ satisfies the
inequalities $2^r<n\leq 2^{r+1}$, which for dimensional reasons
leads to the desired relation $(1+w_2+w_3)^{2^{r+3}}=1$. Now it is
not difficult to see that
$$
g_k=\sum_{k/3\leq i\leq k/2}{i \choose 3i-k}w_2^{3i-k}w_3^{k-2i}.
$$

By an easy computation, $g_{13}=0, g_{12}=w_3^4+w_2^6$ and
$g_{11}=w_2^4w_3$. It turns out that
$w_3^5=w_3g_{12}+w_2^2g_{11}\in J_{n,3}$, but $w_2w_3^4\notin
J_{n,3}$, since it does not belong to the span of
$w_2g_{12}=w_2w_3^4+w_2^7$ and $w_3g_{11}=w_2^4w_3^2$.
Consequently, $\overline{w}_{14}(Y)=\tilde{w}_2\tilde{w}_3^4\neq
0$, and by Corollary~\ref{cor:lepo-1},
$N(\tilde{G}_3(\mathbb{R}^{13})\geq 1+2\cdot 30+2\cdot 14=89$.

\section{Concluding remarks}

\subsection{Embeddings with multiple regularity}

The problem of estimating the invariant $N(M^n)$ was in
\cite{GS-1} (see also \cite{GS-thesis}) incorporated in a more
general question of studying $(k,l)$-regular embeddings. By
definition a smooth map $f : M^n \rightarrow  \mathbb{R}^N$ is
$(k, l)$-regular if for every collection of $k+l$ distinct points
$x_1,\ldots, x_k, y_1,\ldots, y_l$ in $M^n$ and $l$ tangent lines
$L_i \subset T_{y_i}(M^n), \, i = 1,\ldots,l$, the set of points
and lines
$$ f(x_1),\ldots, f(x_k),\, df(L_1),\ldots, df(L_l)$$
is affinely independent.

\medskip
When $l = 0$, the notion of $(k, l)$-regularity reduces to affine
$(k-1)$-regularity in the sense of Handel and Segal
\cite{Han-Seg}. On the other hand, a smooth map is $(0,2)$-regular
if and only if it is totally skew.

\medskip
The existence of a $(k,l)$-regular embedding $f : M^n\rightarrow
\mathbb{R}^N$ implies, essentially by the arguments of
Section~\ref{sec:vector-bundle}, a splitting
$$\varepsilon^N \cong \pi^\ast T(F_{l}(M))\oplus \varepsilon^{k+l-1}\oplus \nu$$
of the trivial $N$-dimensional bundle over the configuration space
$F_{k+l}(M)$ of all ordered collections of $k+l$ distinct points
in $M^n$. By definition $\pi^\ast T(F_l(M))$ is the pull-back of
the tangent bundle $T(F_l(M))$ along the projection map $\pi :
F_{k+l}(M)\rightarrow F_l(M)$ and $\nu$ is a bundle of dimension
$N-(n+1)l-k+1$.

This is a clear indication that the problem of studying
$(k,l)$-regular embeddings is amenable to the methods of
Sections~\ref{sec:characteristic-classes} and we hope to return to
this question in a subsequent publication.

\subsection{Open problems}

In this section we collect some open problems pointing to some of
the most interesting questions about totally skew embeddings of
manifolds.

\begin{prob}{\em Determine the exact values of $N(S^2)$ and
$N(\mathbb{R}P^2)$. More generally what is the exact value of
$N(M^2)$ if $M^2$ is a closed or open surface? According to
\cite{Gho-Tab} the only known result is $N(M)=6$ where $D^2\subset
M \subset \mathbb{R}^2$.}
\end{prob}

\medskip An immersion $\phi : M^n \looparrowright \mathbb{R}^N$ is
called {\em totally skew} if whenever $\phi(x)\neq \phi(y)$ the
affine subspaces $d\phi(T_x(M^n))$ and $d\phi(T_y(M^n))$ are skew.
If $f : M\rightarrow \mathbb{R}^N$ is a totally skew embedding and
if $g: \widetilde{{M}}\rightarrow M$ is a covering map then $\phi
= f\circ g$ is clearly a totally skew immersion. The following
conjecture reflects our impression that in this case a totally
skew immersion can be perturbed to a genuine totally skew
embedding.

\begin{conj}{\em If $M$ is a closed, smooth manifold and $\Gamma$ a
finite group acting freely on $M$ then
$$
N(M)\leq N(M/\Gamma),
$$
in particular $N(S^n)\leq N(\mathbb{R}P^n)$.}
\end{conj}

\medskip
It follows from the splitting (\ref{eqn:splits}) that the
geometric dimension $g$-dim$(\nu)$ of the normal bundle $\nu_1 =
\nu(T(F_2(M)))$ satisfies the inequality
$$
N(M^n)-1 \geq 2n + g\mbox{\rm -dim}(\nu_1).
$$
Similar inequalities hold for manifolds $X, Y$ and $X\times Y$ and
their comparison suggests the possibility of the following general
result.

\begin{conj}{\em For two manifolds $X$ and $Y$, $N(X\times Y)\geq
N(X)+N(Y)-1$.}
\end{conj}

This conjecture, if true, would together with the bound
$N(\mathbb{R}^n) \geq 3n$ for $n$ a power of $2$ (obtained in
\cite{Gho-Tab}), imply the lower bound $N(\mathbb{R}^n)\geq
3n-\alpha (n)+1$.

\medskip
The well-known {\em Immersion Conjecture}, resolved positively by
R.~Cohen \cite{Cohen} in 1985, predicted that any compact smooth
$n$-manifold for $n>1$ can be immersed in
$\mathbb{R}^{2n-\alpha(n)}$, where $\alpha(n)$ is the number of
non-zero digits in the binary representation of $n$. The following
result of Massey, which preceded Cohen's theorem by 15 years,
played an important role by providing strong evidence in favor of
the conjecture.

\begin{theo}{\rm (W.S.~Massey, \cite{Massey})}\label{thm:Massey}
Let $M^n$ be a smooth, compact $n$-dimensional manifold $(n
> 1)$. Then  $\overline{w}_j(M) = 0$
for $j > n -\alpha(n)$.
\end{theo}
Theorem~\ref{thm:Massey} together with our
Corollary~\ref{cor:lepo-1} provides interesting initial evidence
for the following bold conjecture.

\begin{conj}\label{conj:najinteresatnije}
For every $n$-dimensional, compact smooth manifold $M^n$ $(n>1)$,
$$N(M^n)\leq 4n-2\alpha (n)+1.$$
\end{conj}

If correct, Conjecture~\ref{conj:najinteresatnije} would, together
with Proposition~\ref{prop:proj-donja-ocena} and
Theorem~\ref{thm:proj-product}, yield some exact computations of
the invariant $N(M^n)$. For example it would imply
$$N(\mathbb{R}P^2)=7$$
and more generally the following result.
\begin{conj}
\label{wild} Suppose that $n_i=2^{r_i}\; (i=1,...,k)$ and assume
that $r_i\neq r_j$ for $i\neq j$. Let $n=n_1+\cdots +n_k\geq 2$.
Then,
$$N(\mathbb{R}P^{n_1}\times \cdots \times \mathbb{R}P^{n_k})=4n-2\alpha (n)+1.$$
\end{conj}

\vfill \newpage

{\small \DJ{}OR\DJ{}E BARALI\'{C}, Mathematical Institute SASA,
Kneza Mihaila 36, p.p.\ 367, 11001 Belgrade, Serbia

E-mail address: djolebar@nadlanu.com

\bigskip
BRANISLAV PRVULOVI\'{C}, Faculty of Mathematics, University of
Belgrade, Studentski Trg 16, 11000 Belgrade, Serbia

E-mail address: bane@matf.bg.ac.rs

\bigskip
GORDANA STOJANOVI\'{C}, Faculty of Mathematics, University of
Belgrade, Studentski Trg 16, 11000 Belgrade, Serbia

E-mail address: gstojanovic@matf.bg.ac.rs

\bigskip
SINI\v{S}A VRE\'{C}ICA, Faculty of Mathematics, University of
Belgrade, Studentski Trg 16, 11000 Belgrade, Serbia

E-mail address: vrecica@matf.bg.ac.rs

\bigskip
RADE \v{Z}IVALJEVI\'{C}, Mathematical Institute SASA, Kneza
Mihaila 36, p.p.\ 367, 11001 Belgrade, Serbia

E-mail address: rade@mi.sanu.ac.rs}

\end{document}